\documentclass[11pt]{article}
\usepackage{authblk}
\usepackage{hyperref}
\usepackage{amsmath,amssymb, amscd}
\usepackage{amsfonts}
\usepackage{multicol}
\usepackage[latin1]{inputenc}
\usepackage{graphicx}
\usepackage{color}
\usepackage{cleveref}
\usepackage{fancybox}
\usepackage[]{fancyhdr}
\usepackage{enumerate}
\usepackage{makeidx}
\usepackage[all]{xy}
\usepackage{array}
\usepackage{calc}
\usepackage{indentfirst}
\usepackage{url}
\usepackage[shortlabels]{enumitem}

\usepackage{xcolor}
\usepackage{tikz}
\usetikzlibrary{snakes,arrows,shapes}
\newcommand{\dem}{\noindent \textbf{Proof.}\quad}
\newcommand{\fim}{\hfill $\square$ \vspace{0.2cm}}

\newcommand{\x}{\mathbf{v}} % representa um autovetor

\newcommand{\e}{\text{\large{$\mathbf{j}$}}} % representa o vetor com todas as entradas 1.

\newcommand{\ms}{MainSpec} % parte do espectro formada somente pelos autovetores principais.

\newtheorem{teo}{Theorem}%[chapter]

\newtheorem{lema}[teo]{Lemma}

\newtheorem{prop}[teo]{Proposition}

\newtheorem{coro}[teo]{Corollary}

\begin{document}
\title{On main  eigenvalues of certain graphs}

\author[1]{Nair  Abreu}
\author[2,3]{Domingos M. Cardoso}
\author[4]{Francisca A.M. Fran\c{c}a}
\author[5]{Cybele T.M. Vinagre}

\affil[1]{\small
PEP/COPPE, Universidade Federal do Rio de Janeiro, Rio de Janeiro, Brasil. e-mail: nairabreunovoa@gmail.com}
\affil[2]{\small
Center for Research and Development in Mathematics and Applications}
\affil[3]{\small
Department of Mathematics, Universidade de  Aveiro, 3810-193, Aveiro, Portugal. e-mail: dcardoso@ua.pt}
\affil[4]{\small
Instituto de Ci\^{e}ncias Exatas, Universidade Federal Fluminense, Volta Redonda, Brasil. e-mail: francisca@puvr.uff.br}
\affil[5]{\small
Instituto de Matem\'atica e Estat\'{\i}stica, Universidade Federal Fluminense, Niter\'oi, Brasil. e-mail: cybl@vm.uff.br}
\maketitle

\begin{abstract}
An eigenvalue of the adjacency matrix of a graph  is said to be   \emph{main}  if the all-1 vector is not orthogonal
to the associated eigenspace. In this work, we approach the main eigenvalues of some graphs. The
graphs with exactly two main eigenvalues are considered and a relation between those main eigenvalues is presented.
The particular case of harmonic graphs is analyzed and they are characterized in terms of their main eigenvalues
without any restriction on its combinatorial structure.
We give a necessary and sufficient condition for a graph  $G$ to have $-1-\lambda_{\min}$ as an eigenvalue of its
complement, where $\lambda_{\min}$ denotes the least eigenvalue of $G$. Also, we prove that among  connected bipartite
graphs, $K_{r,r}$ is the unique graph for which the index of the complement is equal to $-1-\lambda_{\min}$. Finally,
we characterize all paths and all double stars (trees with diameter three) for which the smallest eigenvalue is non-main.
Main eigenvalues of paths and double stars are identified.
%\keywords{Main eigenvalue \and Harmonic graph \and Path\and Double star.}
 %\subclass{  05C50 \and 05C05 \and 05C38}
\end{abstract}


\noindent \textbf{Classification MSC (2010):} 05C50, 05C05, 05C38.

\section{Introduction and preliminaries}
\label{intro}
In 1970, Cvetkovi\'{c} \cite{cvet1970} introduced the concept of main eigenvalue of a graph, that is, the ones for which
the associated eigenspaces are non-orthogonal to the vector whose entries equal $1$. One year later, the same author \cite{cvet1971}
related the main eigenvalues directly to the number of walks in a graph. It is well known that a graph is regular if and
only if it has only one main eigenvalue \cite{Row2007} but the characterization of graphs with exactly $s >1$ main eigenvalues
is a problem proposed in \cite{cvetk1978} which remains open. The graphs with just two main eigenvalues, have been studied in
several papers \cite{Lepovic2001},\cite{H2002}, \cite{teranishi}, \cite{Hou2006}, \cite{Niki2006}, \cite{Row2007}. Namely, in \cite{Row2007}, a complete survey
on main eigenvalues of a graph (up to 2007) is given and the graphs with exactly two main eigenvalues deserve particular attention.
New results on these graphs appear in more recent publications \cite{Card2008}, \cite{Shi2009}, \cite{hu_li_zhu2009}.\\

In Section~\ref{sec:1}, some properties of the graphs with exactly two main eigenvalues are
studied. Namely, a relation between the two main eigenvalues is deduced and the harmonic graphs are characterized
by its main spectral properties, without any restriction on its combinatorial structure.\\

Em 1971, Cvetkovi\'{c} \cite{cvet1971} started to investigate the relation between the main eigenvalues of a graph
and the eigenvalues of its complement. This subject was also approached in \cite{H2002}, \cite{teranishi} and
\cite{CardPinhe2009}. In \cite{CardPinhe2009}, such relations were investigated seeking
estimates for the maximal size of regular induced subgraphs in the context of convex quadratic programming
applied to graphs. At the end of the paper, three questions are raised. The first wonders about the existence
of a connected graph $G$ whose complement $\overline{G}$ has a main eigenvalue between its first eigenvalue
(\emph{index}) and $-1-\lambda_{\min}(G)$, where $\lambda_{\min}(G)$ denotes the least eigenvalue of $G$. The
second raises the possibility of characterizing graphs $G$ whose spectrum of $\overline{G}$ contains
$-1-\lambda_{\min}(G)$ as a main eigenvalue. The third question approaches the characterization of connected
graphs for which the least eigenvalue is non-main. To answer these questions was the first motivation for the
results of Sections~\ref{sec_3} and \ref{sec_4}.\\

In Section~\ref{sec_3} we answer, in the negative form, to the first question posed
in \cite{CardPinhe2009}. The largest and the second largest eigenvalues of the complement of a graph
are also analyzed and we conclude that $-1-\lambda_{\min}(G)$ belongs to its spectrum if and only if
it coincides with its second largest eigenvalue. We show that, among all connected bipartite graphs,
the  balanced complete bipartite graphs $K_{r,r}$ are those whose respective complements contain
$-1-\lambda_{\min}(G)$ as a main eigenvalue.

In Section~\ref{sec_4} we determine the \emph{main spectrum}  (the set of distinct main eigenvalues)
of a path with $n$ vertices and conclude that the least eigenvalue of such graph is non-main if and
only if $n$ is even. Finally, we conclude that among the trees of diameter three (double stars) only
the balanced ones have the least eigenvalue non-main. On the other hand, the main eigenvalues of an
arbitrary double star are determined and it is shown that their main spectra has cardinality four
when they are not balanced.

Throughout this paper, unless otherwise stated, $G$ denotes a simple graph of order $n$ with edge set $E(G)$ and
vertex set $V(G)=\{ 1, \cdots , n \}$. The edges with end-vertices $i$ and $j$ are simply denoted $ij$ and the
complement of the graph $G$ is denoted $\overline{G}$.
The \label{te1}\textit{adjacency matrix} of $G$, $ \mathbf{A}=[a_{ij}]$, is the  $n \times n$  matrix for which
the entries are $a_{ij}=1$ if \ $ij \in E(G)$, and 0 otherwise. The eigenvalues of $\mathbf{A}$ are also called
the \textit{eigenvalues of} $G$. We write \label{te3}$Spec(G)$ for the multi-set of eigenvalues of $G$. The
characteristic polynomial of $\mathbf{A}$ is called the \emph{characteristic} \emph{polynomial}\emph{ of} $G$.
Unless otherwise stated, the eigenvalues of $G$ are considered in non-increasing order, that is,
\label{te4}$\lambda_{max}=\lambda_1 \ge \lambda_2 \ge \dots \ge \lambda_n = \lambda_{min}$. When necessary,
we write $\mathbf{A}(G)$ instead of $\mathbf{A}$ and $\lambda_{i}(G)$ instead of $\lambda_{i},$ for  $i \in \{1, 2, \dots, n\}$.
The eigenvectors associated to the eigenvalues of $G$ are also called the \textit{eigenvectors of} $G$ and the
eigenspace associated to the eigenvalue $\lambda$ of $G$ is denoted \label{te5}$\varepsilon_{G}(\lambda)$.
 An eigenvector associated to the largest eigenvalue of $G$ is usually called \textit{principal
eigenvector} of $G$.
The all one ${n\times n}$ matrix is denoted \label{te00}$\textbf{J}$ and \label{te0}{\e} denotes
a column of the matrix $\textbf{J}$. An eigenvalue $\lambda$ of $G$ is said to be a \textit{main}
\emph{eigenvalue} if there is an associated eigenvector $\x$ which is not orthogonal to \e. Otherwise,
we say that  $\lambda$ is a \textit{non-main} eigenvalue of $G$.
 Notice that for every graph $G$, its largest eigenvalue $\lambda_1$ is a main eigenvalue.
In particular, when $G$ is $r$-regular with spectrum $\lambda_1=r, \lambda_2, \ldots , \lambda_n$, all
eigenvalues but $\lambda_1$ are non-main \cite{cvet1979}. %\marginpar{Esta parte do texto mudou do final da secção para aqui}

For the basic notions and notation from spectral graph theory not herein defined the  reader is referred to \cite{cvet2010}.
For further mention, we also recall the following consequence of the theorem of Perron-Frobenius (see \cite{cvet1997}).

\begin{teo}\label{teo116}
A graph $G$ is connected if and only if its largest eigenvalue is simple and there exists an associated eigenvector
for which all the coordinates are positive.
\end{teo}

\section{Graphs with exactly two main eigenvalues}
\label{sec:1}
In this section, a relation between the main eigenvalues of a graph with exactly two main eigenvalues
is presented and special attention is given to the harmonic non regular graphs which are characterized using just
their main spectral properties. Furthermore, the bipartite graphs with just two main eigenvalues are also analyzed.
From now on, $\mathbf{d}_G=[d_1, \dots, d_n]^\top$ denotes the degrees vector of  $G$, where $d_i$ is the degree of
the vertex $i \in V(G)$.
\subsection{A relation between the two main eigenvalues}
\label{sec:2}
 \begin{prop}\label{index_relation}
If $G$ is a graph with $m$ edges and exactly two main (distinct) eigenvalues $\lambda_1$ and $\lambda_i$ then
\begin{equation}
\lambda_i = \frac{\sum_{i \in V(G)}{d_i^2}-2m\lambda_1}{2m-n\lambda_1} \label{main_eq_1}\,.
\end{equation}
\end{prop}

\dem
It is immediate that there are scalars $\alpha$ and $\beta$ and orthonormal eigenvectors $\mathbf{u}$
and $\mathbf{v}$ of $\mathbf{A}$ associated to $\lambda_1$ and $\lambda_i$, respectively, such that
$\mathbf{j} = \alpha \mathbf{u} + \beta \mathbf{v}$ and then
$\mathbf{A}\mathbf{j} = \mathbf{d}_G = \alpha \lambda_1 \mathbf{u} + \beta \lambda_i \mathbf{v}$. Then,
it follows
\begin{center}
\begin{tabular}{lcrcrclc}%\hline
$\mathbf{j}^\top\mathbf{j}$    &$=$&$\alpha^2$           &$+$&$\beta^2$           &$=$&$n$  & (i)\\
$\mathbf{j}^\top\mathbf{d}_G$  &$=$&$\lambda_1\alpha^2$  &$+$&$\lambda_i\beta^2$  &$=$&$2m$ & (ii)\\
$\mathbf{d}_G^\top\mathbf{d}_G$&$=$&$\lambda_1^2\alpha^2$&$+$&$\lambda_i^2\beta^2$&$=$&$\sum_{i \in V(G)}{d_i^2}$ & (iii).\\
\end{tabular}
\end{center}
Let us consider two cases (a) $\lambda_i=-\lambda_1$ and (b) $\lambda_i > -\lambda_1$.
\begin{itemize}
\item[(a)] Replacing $\lambda_i$  by $-\lambda_1$ in  (iii), it follows that $\lambda_1^2\left(\alpha^2+\beta^2\right)=\sum_{i \in V(G)}{d_i^2}$.
           Therefore, applying (i), we obtain $\lambda_1^2=\frac{\sum_{i \in V(G)}{d_i^2}}{n}$ which in this case is equivalent to
           \eqref{main_eq_1}.
\item[(b)] Considering the pairs of equations (i)-(ii) and (i)-(iii), we obtain
           $$
           (i)-(ii) \left\{\begin{array}{lcl}
                       \alpha^2& = & \frac{2m-n\lambda_i}{\lambda_1 - \lambda_i}\\
                       \beta^2 & = & \frac{n\lambda_1 - 2m}{\lambda_1-\lambda_i}
                 \end{array}\right. \qquad
           (i)-(iii) \left\{\begin{array}{lcl}
                         \alpha^2& = & \frac{\sum_{j \in V(G)}{d_j^2}-n\lambda^2_i}{\lambda^2_1 - \lambda^2_i}\\
                         \beta^2 & = & \frac{n\lambda^2_1-\sum_{j \in V(G)}{d_j^2}}{\lambda^2_1 - \lambda^2_i}\\
                  \end{array}\right. .
           $$
           \noindent  Therefore, from $2m-n\lambda_i = \frac{\sum_{i \in V(G)}{d_i^2}-n\lambda^2_i}{\lambda_1+\lambda_i}$,
           the equality \eqref{main_eq_1} follows.   \fim
\end{itemize}

Let us call the graphs with the same main eigenvalues \textit{co-main-spectral} graphs.
Despite the relation \eqref{main_eq_1} between the index of $G$ and the other main eigenvalue, in \cite{Card2008}
infinite families of non isomorphic co-main-spectral graphs with exactly two main eigenvalues were presented.
For instance, the bidegreed graphs (that is, graphs where all vertices have one of two possible
degrees) $H^q_k$ obtained from a connected $k$-regular graph $H_k$ of order $p$, after attaching $q \ge 1$
pendent vertices to each vertex of $H_k$ (then the order of $H^q_k$ is $n=(q+1)p$) were considered.
All of these graphs (independently of $p$) have exactly the two main eigenvalues
$\lambda_{1,i}(H^q_k) = \frac{k \pm \sqrt{k^2+4q}}{2}.$ If $k=2$, that is, $H_2$ is the cycle $C_p$,
then $\forall p \ge 3$
$$
\lambda_1(H^q_2)=1+\sqrt{1+q} \; \text{ and } \; \lambda_i(H^q_2)=1-\sqrt{1+q},
$$
are its two main eigenvalues. Notice that for $p=3, 4, \ldots$ we obtain an infinite sequence of
co-main-spectral graphs of increasing order equal to $(q+1)p$ (see \cite[Fig. 2]{Card2008}).
In spite of this, taking into account that $H^q_2$ has $m=n$ edges and
$\sum_{i \in V(H^s_2)}{d_i^2} = p(q+2)^2+pq$ it is easy to check that the equality \eqref{main_eq_1} holds.\\

The next immediate corollary of Proposition~\ref{index_relation} it will be useful for the study of harmonic graphs.

\begin{coro}\label{main_0}
If $\lambda_i=0$, then $\lambda_1 = \frac{\sum_{i \in V(G)}{d_i^2}}{2m}.$
\end{coro}

\subsection{Harmonic graphs}

A graph $H$ is said to be \textit{harmonic} when $\mathbf{d}_{H}$ is an eigenvector associated
to a (necessarily) integer eigenvalue, that is, if  there is a positive integer $\ell$ such that
$\mathbf{A} \mathbf{d}_{H}=\ell \mathbf{d}_{H}$.  It is immediate that every regular graph is harmonic.
The harmonic graphs were introduced in \cite{Grune2002}, \cite{Dress2003}. In \cite{Niki2006}, such a graph without
isolated vertices is called \textit{pseudo-regular graph} and it is defined as being a graph $H$ such that
$\sum_{j \in N_H(i)}{\frac{d_H(j)}{d_H(i)}}$ is constant for every $i \in V(H)$.

The particular case of harmonic trees was studied in \cite{Grune2002}, where the author
consider the trees $\mathcal{T}_{\ell},$ with $\ell \ge 2,$ such that one of its vertices $v$ has degree
$\ell^2 - \ell + 1$, while every neighbor of $v$ has degree $\ell$ and all the remaining vertices have
degree 1. He proved that these are the unique harmonic trees. These trees are among the trees with two
main eigenvalues which have been characterized in \cite{Hou2005} (see also \cite{Hou2006}).

\begin{teo}[\cite{Hou2005}]\label{houarv}
The stars, the balanced double stars and the harmonic trees $\mathcal{T}_{\ell}$, for $\ell \geq 2,$
are the unique trees with exactly two main eigenvalues.
\end{teo}

\begin{lema}\label{bipertite_harmonic}
If $H$ is a bipartite harmonic graph with at least one edge and largest eigenvalue $\lambda_1$,
then $-\lambda_1$ is a non-main eigenvalue of $H$.
\end{lema}

\dem
Let us consider that the bipartite harmonic graph $H$ has $q \ge 1 $ connected nontrivial components
$H_1, \ldots, H_q$ and $H_{q+1}, \ldots, H_{q+p}$ trivial components, with $p \ge 0$ . It is immediate
that each $H_k$ of order $n_k$ is a connected bipartite harmonic subgraph with the same largest eigenvalue
$\lambda_1$ and the same simple least eigenvalue $-\lambda_1$. Assuming that $V(H_k)$ admits the bipartition
$S_k$ and $T_k$ such that each edge of $H_k$ has one end-vertex in $S_k$ and the other in $T_k$, then the vectors
$\left(\begin{array}{c}
                         d_{S_k}\\
                         d_{T_k}\\
                   \end{array}\right)       \text{ and }    \left(\begin{array}{c}
                                                                         -d_{S_k}\\
                                                                          d_{T_k}\\
                                                                 \end{array}\right),
$
where $d_{S_k}$ and $d_{T_k}$ denote the subvectors of degrees of the vertices in $S_k$ and $T_k$,
are the principal eigenvector and the eigenvector associated to $-\lambda_1$, respectively, of $H_k$,
for $k=1, \ldots, q$. The vectors
\begin{eqnarray}
\hat{u}_k^T &=& (\begin{array}{ccccccccccc} 0, & \dots,& 0, &\;d^T_{S_k},& d^T_{T_k}, & 0, & \dots, & 0, & 0, & \dots,  & 0\end{array})\label{vector_u} \\
\hat{v}_k^T &=& (\begin{array}{ccccccccccc} 0, & \dots,& 0, & -d^T_{S_k},& d^T_{T_k}, & 0, & \dots, & 0, & 0, & \dots,  & 0\end{array})\label{vector_v}
\end{eqnarray}
where the last $p$ zero coordinates correspond to the $p$ trivial components, the $k-1$ zeros on the left
of $d^T_{S_k}$ and the $q-k$ zeros on the right of $d_{T_k}$ correspond to the vertices in the components
$H_1, \dots, H_{k-1}, H_{k+1}, \dots, H_{q}$, respectively, are also a principal eigenvector and the
eigenvector associated to $-\lambda_1$, respectively, of $H$. Therefore, since for each component $H_k$
the sum of the degrees of the vertices in $S_k$ is equal to the sum of the degrees of the vertices in
$T_k$, it follows that the vectors in \eqref{vector_v} are all orthogonal to the all one vector. Since every vector of
the eigenspace associated to $-\lambda_1$ is a linear combination of those vectors in \eqref{vector_v},
$-\lambda_1$ is non-main. \fim

Nikiforov in \cite[Th. 8]{Niki2006} proved that every main eigenvalue of an harmonic graph $H$ belongs
to the set $\{-\lambda_1, 0, \lambda_1\}$. It is also stated in \cite[Th. 8]{Niki2006} that if $H$ is
a graph without a bipartite component such that all  main eigenvalues are in
$\{-\lambda_1, 0, \lambda_1\}$, then it is harmonic. A similar result for connected graphs is obtained
in \cite[Pr. 3.3]{Row2007}, using a different approach. The next proposition gives a spectral
characterization of harmonic graphs without any restriction regarding theirs combinatorial structure.

\begin{prop}\label{harmprincipal}
A graph $H$ is harmonic if and only if every main eigenvalue of $H$ belongs to the set $\{0, \lambda_1\}$.
\end{prop}

\dem
If $H$ is harmonic, as direct consequence of Lemma~\ref{bipertite_harmonic} and Theorem 8 in \cite{Niki2006},
it follows that every main of its eigenvalues are in $\{0, \lambda_1\}$.
Converse\-ly, let us consider that the main eigenvalues of $H$ are in $\{0, \lambda_1\}$. If $H$ has only one
main eigenvalue, then $H$ is regular and the result follows. Other\-wise, assuming that
$\{\mathbf{v}_1, \ldots, \mathbf{v}_p\}$ is a basis for $\varepsilon_H(0)$ and
$\{\mathbf{u}_1, \ldots, \mathbf{u}_q\}$ is a basis for $\varepsilon_H(\lambda_1)$, it follows that
$\mathbf{j}=\sum_{i=1}^p{\alpha_i \mathbf{u}_j}+\sum_{j=1}^q{\beta_j \mathbf{v}_j}$ for some scalars
$\alpha_1, \ldots, \alpha_p, \beta_1, \ldots, \beta_q$ and
$\mathbf{d}_H = A_H \mathbf{j} = \lambda_1 \sum_{i=1}^p{\alpha_i \mathbf{u}_j}.$ Therefore,
$\mathbf{d}_H \in \varepsilon_H(\lambda_1)$. \fim

 The following proposition gives an alternative characterization of harmonic graphs.

\begin{prop}
A graph $G$  with $m$ edges is harmonic if and only if $\lambda_1 = \frac{\sum_{i \in V(G)}{d_i^2}}{2m}$
and it has no more than two main eigenvalues.
\end{prop}

\dem
Suppose that $G$ is harmonic. By Proposition~\ref{harmprincipal}, all its main eigenvalues are in $\{0, \lambda_1\}$
and we have two cases: (i) $G$ is regular, with degree say $k$, and then $\lambda_1=k$ is the unique main eigenvalue
or (ii) $G$ is non regular and then it has two main eigenvalues.
\begin{itemize}
\item[(i)] $\lambda_1=k=\frac{nk^2}{nk}=\frac{\sum_{i \in V(G)}{d_i^2}}{2m}$.
\item[(ii)] By Corollary~\ref{main_0}, it follows that $\lambda_1 = \frac{\sum_{i \in V(G)}{d_i^2}}{2m}$.
\end{itemize}
Conversely, assume that $\lambda_1 = \frac{\sum_{i \in V(G)}{d_i^2}}{2m}$ and $G$ has no more than two main
eigenvalues. If $G$ is regular then the conclusion is immediate. Else, by Proposition~\ref{index_relation}, the
main eigenvalues of $G$, $\lambda_i(G)$ and $\lambda_1$, are related by the equality \eqref{main_eq_1}.
Replacing $\lambda_1$ in \eqref{main_eq_1} by $\frac{\sum_{i \in V(G)}{d_i^2}}{2m}$ it follows that the
main eigenvalues of $G$ are in $\{0,\lambda_1\}$. Therefore,  by Proposition~\ref{harmprincipal}, the result
follows. \fim

\section{The largest and the second largest eigenvalues of the complement of a graph}\label{sec_3}
From now on, we consider the all distinct eigenvalues $\mu_1, \dots , \mu_s$,  $1 \le s \le n$, of the graph
$G$ having the respective associated eigenspace %$\varepsilon_{G}(\mu_i)$,
not orthogonal to the vector $\e$ as the \emph{main}
eigenvalues of $G$ and the remaining distinct eigenvalues $\mu_{s+1}, \cdots , \mu_p$, $s+1 \le p \le n$, as
the \emph{non-main} eigenvalues.
The set of distinct main eigenvalues of $G$ is herein called the \emph{main spectrum} of $G$ and it is denoted
$\ms(G)$. Therefore,
$Spec(G) = \{\mu^{[q_1]}_1, \dots , \mu^{[q_s]}_s, \mu^{[q_{s+1}]}_{s+1}, \dots , \mu^{[q_p]}_s\}$,
where $\mu^{[q_j]}_j$ means that the eigenvalue $\mu_j$ has multiplicity $q_j$.

Before to proceed, it is worth to recall the following theorem.

\begin{teo}[\cite{cvet1971}]\label{cvet71}  $\ms(G)$ and $\ms(\overline{G})$ have the same number of elements.
Furthermore, if $\lambda \in \ms(G)$ and   $\overline{\lambda}   \in \ms(\overline{G})$, then
$\lambda +\overline{\lambda} \neq -1$ .
\end{teo}

Taking into account this theorem and the definition of main/non-main eigenvalue it is immediate to obtain the
basic results stated in the next proposition partially proved in \cite{H2002}.

\begin{prop}\label{Hteo33}
Consider a graph $G$ and $\lambda \in Spec(G)$. Then the following assertions are equivalent:
\begin{enumerate}%[(i)]
\item the eigenvalue $\lambda $ is non-main or it is main with multiplicity greater than $1$; \label{teo231}
\item there is some eigenvector $\x$ of $G$ associated to $\lambda$ such that $\e^\top\x=0$;\label{teo232}
\item the scalar $-1-\lambda$ belongs to $Spec(\overline{G})$.\label{teo233}
\end{enumerate}
\end{prop}

As direct consequence  of this proposition, we may note that a necessary and sufficient condition for a simple
eigenvalue $\lambda$ of a graph $G$ to be non-main is $-1-\lambda$ to be an eigenvalue of $\overline{G}$
(see \cite{H2002}).

Furthermore, we also may conclude the following corollary of Proposition~\ref{Hteo33}.

\begin{coro} \label{coronovo}
If $-1-\lambda(G)$ is a simple eigenvalue of $\overline G$ then it is non-main.
\end{coro}

\dem
It follows from Proposition \ref{Hteo33}, in view of the known relation
$\mathbf{A}(\overline G)=\mathbf{J}-\mathbf{I}_n-\mathbf{A}( G)$. \fim

Now, it is worth to recall the following consequence of Weyl's inequalities which proof can be found in \cite{cvet1997}:
\begin{equation}\label{-1-ln}
{\lambda}_{2}({\overline{G}})\leq -1-\lambda_{n}(G) \le {\lambda}_{1}({\overline{G}}).
\end{equation}

The relations  (\ref{-1-ln})  furnish a (negative) answer to the  question raised in \cite{CardPinhe2009}
about the existence  of a graph $G$ for which the complement $\overline{G}$ has an eigenvalue less than
its index and greater than $-1-\lambda_n(G)$.

\begin{prop}\label{res1}
If $G$ is a graph of order $n$, then $\overline{G}$ has no eigenvalue belonging to the  open interval
$(-1-\lambda_{n}(G), \   \lambda_1({\overline{G}}))$.
%$\overline{\lambda}$ such that
%$-1-\lambda_{n}(G)< \overline{\lambda} < \lambda_1({\overline{G}}).$
\end{prop}

The inequalities \eqref{-1-ln}  motivate us to consider graphs $G$ for which $-1-\lambda_n(G)$ is an eigenvalue
of $\overline G$. We have two cases: (a) $\lambda_1(\overline G)=-1-\lambda_n(G)$ and (b)
$\lambda_2(\overline G)=-1-\lambda_n(G)$.

In the case (a), we have that $-1-\lambda_n(G)=\lambda_1(\overline{G})$ is a main eigenvalue of
$\overline{G}$. Therefore,  Theorem~\ref{cvet71} guarantees that $\lambda_n(G)$ is a
non-main eigenvalue, since $-1= \lambda_n(G)+ (-1-\lambda_n(G)) $. In fact, regarding the equality (a), we may
establish the following proposition.

\begin{prop}\label{lema43}
Let $G$ be a graph of order $n$. Then $\lambda_1({\overline{G}})=-1-\lambda_n(G)$ if and only if $\lambda_n(G)$
is non-main and the multiplicity of ${\lambda}_{1}({\overline{G}})$ is greater than one.
\end{prop}

\dem
If $\lambda_1(\overline{G})=-1-\lambda_n(G)$,  $\lambda_n(G)$ is non-main and from
Proposition~\ref{Hteo33}, $-1-\lambda_n$ has an eigenvector $\x_1$ such that $\x_1 \in \varepsilon_G(\lambda_n)$
and ${\e}^\top \x_1=0$. On the other hand (by Perron-Frobenius theorem), there is an eigenvector $\x_2$,
associated to $\lambda_{1}(\overline{G})$, with nonnegative entries and then $\e^\top \x_2 \ne 0$. Therefore,
$\x_1$ and $\x_2$ are linearly independent. This implies that the multiplicity of $\lambda_1(\overline{G})$ is
greater than $1$. Conversely, if $\lambda_n(G)$ is non-main then $-1-\lambda_n(G) \in Spec(\overline{G})$ by
Proposition \ref{Hteo33}. Since the multiplicity of $\lambda_1(\overline{G})$ is greater than one, from  (\ref{-1-ln})
the result follows.   \fim

According to  Theorem \ref{teo116} and Proposition~\ref{lema43},  when $\lambda_{1}(\overline{G})=-1-\lambda_{n}(G)$
it follows that $\lambda_{n}(G)$ is non-main and $\overline{G}$ is disconnected. On the other hand, for the case (b)
we have:

\begin{prop}\label{prop13}
Let $G$ be a graph of order $n$. Then $\lambda_2(\overline{G})=-1-\lambda_n(G) < \lambda_1({\overline{G}})$
if and only if $\lambda_n(G)$ is main with multiplicity greater than one or it is non-main and
${\lambda}_{1}({\overline{G}})$ is simple.
\end{prop}

The inequalities in (\ref{-1-ln}) and Propositions \ref{lema43} and \ref{prop13}  allow us to conclude that
$-1-\lambda_{n}(G) \in Spec(\overline{G})$ if and only if ${\lambda}_2(\overline{G})=-1-\lambda_n(G)$.

From Corollary \ref{coronovo}, for an arbitrary  graph $G$  of order $n$  such that $\lambda_{n}(G)$ is a simple
eigenvalue, we have that  $\lambda_{n}(G)$ is non-main if and only if $-1-\lambda_{n}(G)$ is an eigenvalue of
$\overline{G}$. Since the least eigenvalue of a connected bipartite graph is simple, for these graphs we may
conclude the following:
\begin{itemize}
\item[(a)] For a connected bipartite graph $G$, $\lambda_1(\overline{G})=-1-\lambda_n(G)$ if and only if $\lambda_n(G)$ is non-main and
           ${\lambda}_{1}({\overline{G}})$ has multiplicity greater than one.

\item[(b)] If $G$ is connected and bipartite then $\lambda_2(\overline{G})=-1-\lambda_n(G)< \lambda_1({\overline{G}})$ if and only if $\lambda_n(G)$
           is non-main and ${\lambda}_{1}({\overline{G}})$ is simple.
\end{itemize}

The next result gives a combinatorial characterization of bipartite graphs $G$ of order $n$ for which ${\lambda}_{1}({\overline{G}})=-1-\lambda_{n}(G)$.

\begin{teo}\label{res2}
Let $G$ be a bipartite graph of order $n$. Then $\lambda_1({\overline{G}})=-1-\lambda_n(G)$ if and
only if  $G$ is complete (bipartite) and balanced.
\end{teo}

\dem
Let  us consider a bipartite graph $G$ with vertex set  $V=V_1\dot{\cup}V_2$, where
$|V_1|=r$ and $|V_2|=s$. If $\lambda_1({\overline{G}})=-1-\lambda_n(G)$ then (a) above implies
$\overline{G}$ is disconnected, and thus $\overline{G}=K_r \dot{\cup}K_s$. Since $\lambda_1(G)$
is a multiple eigenvalue then $r=s$.
 Conversely, if  $G=K_{s,s}$, for some positive integer $s$, then $\overline{G}=\overline{K_{s,s}}$
is a disconnected graph with two components which are complete graphs with $s$ vertices. It follows that
$\lambda_n(G)=-s$ and $\lambda_1({\overline{G}})=s-1$ and therefore,    $\lambda_1(\overline{G})=-\lambda_n(G)-1$. \fim

\section{Main spectra of paths and double stars.}\label{sec_4}

Concerning the third question of \cite{CardPinhe2009}, we may note that among the connected graphs for which
the least eigenvalue is non-main we can count the harmonic graphs (see the Proposition~\ref{harmprincipal})
which includes the regular graphs. In this section, the graphs with non-main least eigenvalue of two
families of trees are characterized. We start by determining the paths with non-main least eigenvalue. For
sake of completeness, we determine the main spectrum of an arbitrary path.

It is worth to recall the following lemma which can be found in \cite {cvet1979} (the eigenvectors are
described in \cite{Lee1993}).

\begin{lema}[\cite {cvet1979},\cite{Lee1993}] \label{lee1993}
Let $\mathcal{P}_n$ be the path on $n$ vertices. Then its eigenvalues are simple and given by
$\lambda_j(\mathcal{P}_n)=2\cos\left(\dfrac{j\pi}{n+1}\right)$, $1 \leq j \leq n$.
Each of these eigenvalues $\lambda_j$ has an associated eigenvector with entries
$\x^{(j)}_{i}={\sin}\left(i\dfrac{j\pi}{n+1}\right)$, for $i\in\{1, 2, \dots, n\}$.
\end{lema}

\begin{teo}\label{qcap} For $n \ge 2$ and
$1 \leq j \leq n$, $\lambda_j$ is a non-main eigenvalue of the path $\mathcal{P}_{n}$ if and only if $j$
is even. In particular, the least eigenvalue of  $\mathcal{P}_{n}$ is non-main if and only if $n$ is even.
\end{teo}

\dem
Let us fix $j$, $1 \leq j \leq n$. For the $\lambda_j$-eigenvector
$\x^{(j)}=(\x^{(j)}_1, \ldots , \x^{(j)}_n)^{\top}$ we have $\lambda \x^{(j)}_i =  \sum_{t\sim i} \x^{(j)}_t$,
whence $\lambda\sum_{i}\x^{(j)}_i = \sum_{i}d_i\x^{(j)}_i = 2 \sum_{i} \x^{(j)}_i -\x^{(j)}_1 -\x^{(j)}_n$.
From Lemma~\ref{lee1993}, $\lambda_j \neq 2$ and then $\sum_{i} \x^{(j)}_{i} = 0$ if and only if
$\x^{(j)}_1 + \x^{(j)}_n = 0$. Since
$\x^{(j)}_1 + \x^{(j)}_n = 2\sin\left(\frac{j\pi}{2}\right)\cos\left(\frac{(n-1)j\pi}{2(n+1)}\right)$,
we may verify that $\lambda_j$ is a non-main eigenvalue if and only if $ j$ is even. In fact,
$\cos\left(\frac{(n-1)j\pi}{2(n+1)}\right)=0$ if and only if $ \frac{(n-1)j\pi}{2(n+1)}=\frac{\pi}{2}+k\pi$,
for $k \in \mathbb{N}$.  From this, we have that $1 + 2k <j \leq n$. Also, it holds  that
$n = \frac{j+(1+2k)}{j-(1+2k)}$,  which implies  $n < \frac{2n}{j-(1+2k)}$ and then,  $j < 2k+3$.
Thus $1+2k < j < 2k+3$, that is, $j$ is even. The another case is straightforward. \fim

\begin{coro}\label{res7}
The path $\mathcal{P}_{n}$ on $n$ vertices has $\lceil \frac{n}{2} \rceil$ main eigenvalues, where
$\lceil x \rceil$ denotes the least integer no less than $x$.
\end{coro}

The following result characterizes the semi-regular bipartite graphs in terms of theirs main eigenvalues.

\begin{teo}[\cite{Row2007}]\label{lnmain}
A non-trivial connected graph $G$ is semi-regular bipartite if and only if its main eigenvalues are only
$\lambda_{1}(G)$ and $-\lambda_{1}(G)$.
\end{teo}

Combining Theorem~\ref{lnmain} with  Proposition~\ref{index_relation}, it follows that when $G$
is a connected semi-regular bipartite graph of order $n$, $\lambda_{1}(G)=\sqrt{\frac{\sum_{i \in V(G)}{d_i^2}}{n}}$.
This is a known result obtained in \cite{Hofmeister83} where it was stated that if a graph $H$ has order $n$,
then $\sum_{i \in V(H)}{d_i^2} \le \lambda^2_{1}(H)n$ and the equality holds if and only if $H$ is a semi-regular
bipartite graph.

Since diameter-2 trees (the stars) are connected semi-regular bipartite graphs, there
exists no diameter-2 tree with  non-main least eigenvalue. Regarding diameter-3 trees, it should be noted
that these trees (the double stars) are not semi-regular bipartite graphs and then, combining  Theorems
\ref{lnmain} and \ref{houarv}, we may conclude that the least eigenvalue of a balanced double star is
non-main. On the other hand,  we claim that   there are no non-balanced double stars with least
eigenvalue non-main.

%%%%%%%%%%%%%%%%%%%%%%%%%%%%%%

In order to prove our assertion, we first remember that a \emph{walk} in $G$ is a sequence
$i_0, i_1, \cdots i_r$ of vertices in $G$ such that $i_t$ is adjacent to $i_{t+1}$, for $0 \leq t \leq r-1$.
The \emph{length} of a walk is its number of edges. For a square matrix $\mathbf{B}$, the  \emph{walk-matrix}
of $\mathbf{B}$ is given by $\mathbf{W(B)} = [\e \  \mathbf{B}\e \   \mathbf{B}^2\e \cdots \  \mathbf{B}^{n-1}\e]$.
In the particular case  $\mathbf{B}=\mathbf{A}$,  the adjacency matrix  of $G$,
$\mathbf{W(A)} = [w_{ij} ]$  is such that $w_{ij}$ gives the number of walks in $G$ of length $j$ starting at vertex $i$,
$1 \leq i \leq n$ and $1 \leq j \leq n-1$, and then it is called the \emph{walk-matrix of} $G$.
We also recall that a partition $\pi$ of the vertex set $V(G)$ of the graph $G$ is  \emph{equitable} when,
given two cells $V_i$ and $V_j$ of $\pi$, there is a constant $m_{ij}$ such each vertex $v \in V_j$ has
exactly $m_{ij}$ neighbors in $V_j$.
The matrix $\mathbf{M}=[m_{ij}]$ is called the \emph{divisor of} $G$ with respect to $\pi$. It is known
(Theorem 3 of \cite{cvetk1978}), that the main eigenvalues of $G$ are eigenvalues of  $\mathbf{M}$.
A fundamental result on the number of eigenvalues of a graph is the following

\begin{teo}[\cite{H2002}]
The rank of the walk-matrix of  $G$ is equal to the number of its main eigenvalues.
\end{teo}

It was recently proved  (\cite{Lulu2015}, Lemma 2.4) that the number of main eigenvalues of $G$ is equal to
the rank of the walk-matrix $\mathbf{W(M)}$  of $\mathbf M$.

\begin{teo}\label{diam3_e_aut_p}
Let $T$ be a double star with $n$ vertices. Then its least eigenvalue is non-main if and only if $T$
is balanced.
\end{teo}

\noindent \dem Let $T=T(k,s)$ be a double star of order $n=k+s+2$ whose vertices are labeled as in Figura~\ref{t(k,s)}.
 \begin{figure}[h]
	\centering
\begin{tikzpicture}[>=latex',join=bevel,scale=0.8]
\tikzstyle{selected edge} = [draw,line width=1pt,->,red!30]
\node (1) at (-30bp,70bp) [label=-90:{\small{$\cdots$}}] {};
\node (4) at (21bp,40bp) [label=-90:{\small{$\cdots$}}] {};
\node (2) at (-30bp,90bp) [draw,circle,inner sep=1.5pt,fill=black!100,label=90:{\footnotesize{$1$}}] {};
\node (3) at (20bp,60bp) [draw,circle,inner sep=1.5pt,fill=black!100,label=90:{\footnotesize{$2$}}] {};
\node (5) at (-45bp,60bp) [draw,circle,inner sep=1.5pt,fill=black!100,label=-90:{\footnotesize{$3$}}] {};
\node (6) at (-15bp,60bp) [draw,circle,inner sep=1.5pt,fill=black!100,label=-90:{\footnotesize{$k+2$}}] {};
\node (7) at (-3bp,30bp) [draw,circle,inner sep=1.5pt,fill=black!100,label=-90:{\footnotesize{$k+3$}}] {};
\node (8) at (45bp,30bp) [draw,circle,inner sep=1.5pt,fill=black!100,label=-90:{\footnotesize{$k+s+2$}}] {};
%\node (9) at (24bp,0bp) [label=-90:{\small{Figura 2: $T=T(k,s)$ tem $n=k+s+2$ vértices}}] {};
 \draw [] (2) -- node {} (3);
 \draw [] (2) -- node {} (5);
 \draw [] (2) -- node {} (6);
 \draw [] (3) -- node {} (7);
 \draw [] (3) -- node {} (8);
\end{tikzpicture}
	\caption{Double star $T(k,s)$.}\label{t(k,s)}
\end{figure}
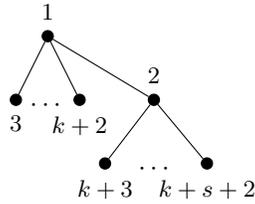
Let us consider $V_1=\{1\}$, $V_2=\{2\}$, $V_3=\{3, \ldots , k+2\}$ and $V_4=\{k+3, \ldots, k+s+2 \}$.
Then $V_1 \dot\cup V_2 \dot\cup V_3 \dot\cup V_4$ is an equitable partition of $V(T)$ with associated divisor
$$
\mathbf{M}=[m_{ij}] =\left[
    \begin{array}{cccc}
      0 & 1 & k & 0 \\
      1 & 0 & 0 & s \\
      1 & 0 & 0 & 0 \\
      0 & 1 & 0 & 0 \\
    \end{array}
  \right]\,,
$$
for which the  walk-matrix is
$$
\mathbf {W(M)} =
  \left[
    \begin{array}{cccc}
      1 & k+1 & k+s+1 & s+k^2+2k+1 \\
      1 & s+1 & k+s+1 & s^2+2s+k+1 \\
      1 & 1 & k+1 & k+s+1 \\
      1 & 1 & s+1 & k+s+1 \\
    \end{array}
  \right].
$$
It can be verified that $\det \mathbf {W(M)}=-ks(s-k)^2$,
which is equal to zero  if and only if $s=k$. Since $\mathbf M$ has characteristic polynomial $q(x)=x^4-(k+s+1)x^2+ks$
and, according to \cite{Vinagre2009}, the characteristic polynomial of $T$ is $p(x)=x^{s-1}x^{k-1}(x^{4}-x^{2}(k+s+1)+ks)$,
we conclude that in case $k \neq s$ the four non-zero eigenvalues of the graph $T=T(k,s)$ are main. In particular,
it follows that $\lambda_{n}$ is a  main eigenvalue (clearly, the others are $\lambda_{1}$, $\lambda_{2}$ and
$\lambda_{n-1}$). Considering that the case $k=r$ is already known, the assertion is proved. \fim

By combining Theorem \ref{cvet71}, Proposition \ref{prop13} and Theorems \ref{qcap} and \ref{diam3_e_aut_p} we may conclude immediately the next corollary.

\begin{coro} If the graph $G$ is a path   (respectively, a balanced double star) on $n$ vertices then its complement $\overline G$ has $\lceil\frac{n}{2}\rceil$ (resp., two) main eigenvalues and the second largest eigenvalue of $\overline G$ is equal to $-1-\lambda_n(G)$.
\end{coro}

\section*{Acknowledgements}
The research of Domingos M. Cardoso is partially supported by the Portuguese Foundation for Science
and Technology (\textquotedblleft FCT-Funda\c c\~ao para a Ci\^encia e a Tecnologia \textquotedblright),
through the CIDMA - Center for Research and Development in Mathematics and Applications, within project
UID/MAT/ 04106/2013. This author also thanks the support of Project Universal CNPq 442241/2014 e Bolsa
PQ 1A CNPq, 304177/2013-0 and the hospitality of PEP/COPPE/UFRJ where this paper was started.

\bibliographystyle{plain}
\bibliography{bibliografia}

\end{document}